\documentclass{amsart}
\DeclareMathSymbol{\twoheadrightarrow} 
{\mathrel}{AMSa}{"10}

\def\Q{{\mathbf Q}}
\def\Z{{\mathbf Z}}
\def\C{{\mathbf C}}

\def\F{{\mathbf F}}

\def\r{{r}}
\def\n{{n}}
\def\A{{A}}

\def\Gal{\mathrm{Gal}}
\def\End{\mathrm{End}}
\def\Aut{\mathrm{Aut}}

\def\GL{\mathrm{GL}}
\def\SL{\mathrm{SL}}

\def\M{\mathrm{M}}

\def\dim{\mathrm{dim}}

\def\invlim{{\displaystyle{\lim_{\leftarrow}}}}

\def\bmu{\boldsymbol \mu}
\newtheorem{thm}{Theorem}[section]
\newtheorem{lem}[thm]{Lemma}
\newtheorem{cor}[thm]{Corollary}
\newtheorem{prop}[thm]{Proposition}
\theoremstyle{definition}
\newtheorem{defn}[thm]{Definition}
\newtheorem{ex}[thm]{Example}
\newtheorem{rem}[thm]{Remark}
\newtheorem{rems}[thm]{Remarks}
\title[Polarizations on abelian varieties]
{Polarizations on abelian varieties}

\author[A.\ Silverberg]{A.\ Silverberg}
\address{Mathematics Department, Ohio State
University,
Columbus, Ohio, 
USA}
\email{silver\char`\@math.ohio-state.edu}

\author[Yu.\ G.\ Zarhin]{Yu.\ G.\ Zarhin}
\address{Mathematics Department, Pennsylvania State
University,
University Park, PA, 
USA\newline\indent
Institute for Mathematical Problems in Biology,
Russian
Academy of Sciences, Push\-chi\-no, Moscow Region,
Russia}
\email{zarhin\char`\@math.psu.edu}

\begin{document}

\maketitle

\section{Introduction}

If $F$ is a field and $\ell$ is a prime number,
let $F^{(\ell)}$ denote the extension of
$F$ obtained by adjoining all $\ell$-power 
roots of unity.
In \S\ref{thm1sect} we prove a result 
(Theorem \ref{eulerthm}) which implies:

\begin{thm}
\label{thm1} 
Suppose $F$ is a number field and $\r$ is a positive 
integer. Then there exists an abelian variety $A$ 
over $F$ with the following properties:
\begin{enumerate}
\item[(a)] every $F$-polarization on
every abelian variety $F$-isogenous
to $A$ has degree divisible by $\r$;
\item[(b)] if $\ell$ is a prime and
$\ell^{a}$ divides $\r$, 
then every $F^{(\ell)}$-polarization
on every abelian variety $F^{(\ell)}$-isogenous
to $A$ has degree divisible by $\ell^{2a}$.
\end{enumerate}
\end{thm}

We deduce Theorems \ref{eulerthm} 
and \ref{thm1} from
Theorem \ref{avprop}, which gives a general 
method for obtaining isogeny 
classes of abelian varieties all of whose 
polarizations have degree divisible by a given 
prime power.
The proof uses some results on elliptic curves
(see Corollary \ref{manyell}) that 
may be of independent interest.

It is known that every isogeny class over an
algebraically closed field contains
a principally polarized abelian variety
(\cite{Mumford}, Cor.~1 to Thm.~4 in \S 23).

Howe (\cite{howetrans}; see also \cite{howejalg})
gave examples of isogeny classes of abelian varieties
over finite fields with no principal polarizations 
(but not with the degrees of all the
polarizations divisible by a given non-zero integer,
as in Theorem \ref{thm1} above).

In \cite{sympl} we obtained,
for all odd primes $\ell$, 
isogeny classes of abelian varieties in positive
characteristic, all of whose 
polarizations have degree divisible by $\ell^2$.
(We gave results in the more general context of
invertible sheaves; see also Theorems \ref{eulerthm}
and \ref{avprop} below.)
Our results gave the first examples for which all the
polarizations of the abelian varieties in an isogeny
class have degree divisible by a given prime.

Inspired by our results, Howe \cite{Howe} recently
obtained, for all odd primes $\ell$, 
examples of isogeny classes of abelian varieties 
over fields of arbitrary characteristic 
(including over number fields),
all of whose polarizations have degree divisible by
$\ell^2$.
The results of this paper were inspired
by an early version (\S 2 of the current version)
of \cite{Howe}. In particular, we employ Howe's idea
of considering the Jordan-H\"older factors of an
appropriate module.

Zarhin was partially supported by
EPSRC grant GR/M98135.
Silverberg would like to thank NSF, NSA, MSRI, 
the Mathematics Institute of the 
University of Erlangen, 
and the Alexander von Humboldt Foundation.

\section{Notation}
\label{prelim}

If $T$ is a $\Z_\ell$-lattice and 
$e : T \times T \to \Z_\ell$ is a non-degenerate
pairing, let 
$$T^* = \{y \in T\otimes_{\Z_\ell}\Q_\ell \mid 
e(T,y) \subseteq \Z_{\ell}\},$$ 
the dual lattice of $T$ with respect to $e$.

Suppose $F$ is a field.
Let $F^s$ be a separable closure of $F$ and
let $G_F=\Gal(F^s/F)$.
If $\A$ is an abelian variety over $F$ and $n$ is a
positive
integer, let $\A_n$ denote the kernel of
multiplication by $n$ in $\A(F^s)$.
Let $T_\ell(\A)= \invlim \A_{\ell^m}$,
let $V_\ell(\A) = T_\ell(\A) \otimes_{\Z_\ell} \Q_\ell$,
let $\bar{\A}=\A\times_F{F^s}$,
let $\A^*$ denote the dual of $\A$,
and let
$$\rho_{A,n} : G_F \to \Aut(\A_n) \cong \GL_{2d}(\Z/n\Z)$$
denote the mod $n$ representation (where $d=\dim(\A)$).
Let $$F^{(\ell)}=F(\bmu_{\ell^{\infty}}) \quad \text{and}
\quad F_\A^{(\ell)}=F(\A_{\ell^{\infty}}).$$

\section{Algebra}
\label{lasect}

\begin{lem}
\label{fulllem}
Suppose $N$ is a group, $S$ is a $\Z_\ell[N]$-module
which is free of finite rank $m$ over $\Z_\ell$,
and $W=S\otimes_{\Z_\ell}\Q_\ell$.
Then the following are equivalent:
\begin{enumerate}
\item[(i)] the natural homomorphism
$\Z_\ell[N] \to \End_{\Z_\ell}(S)$ is surjective;
\item[(ii)]  the $N$-module $S/\ell S$ is absolutely
simple;
\item[(iii)] the image of the natural homomorphism
$\Z_\ell[N] \to \End_{\Q_\ell}(W)$
is isomorphic to the matrix ring $\M_{m}(\Z_\ell)$.
\end{enumerate}
If these equivalent conditions hold, then:
\begin{enumerate}
\item[(a)]
the $N$-module $W$ is absolutely simple;
\item[(b)]
every $N$-stable $\Z_{\ell}$-lattice in $W$
is of the form $\ell^i S$ for some integer $i$;
\item[(c)]
if
$f: S \times S \to \Z_{\ell}$
is an $N$-invariant pairing such that
$f(S,S)=\Z_{\ell}$ then $f$ is perfect.
\end{enumerate}
\end{lem}

\begin{proof}
Assume (ii).
By Theorem 9.2 
of \cite{Isaacs}, the natural homomorphism
 $\F_\ell[N] \to \End(S/\ell S)$ is surjective.
Now (i) follows from Nakayama's Lemma.

Clearly, (i) implies (iii).

Assume (iii).
Then the natural homomorphism  $\F_\ell[N] \to
\End(S/\ell S)$ 
factors through a surjective homomorphism
 $\F_\ell[N] \to \M_m(\F_\ell)$, i.e., $S/\ell S$
is an $\M_m(\F_\ell)$-module. Since $\M_m(\F_\ell)$ is
a simple 
algebra and
$m=\dim_{\F_{\ell}}(S/\ell S)$, the 
map $\M_m(\F_\ell) \to \End(S/\ell S)$ 
is injective and therefore surjective. 
Thus the composition
$\F_\ell[N] \to\M_m(\F_\ell) \to \End(S/\ell S)$ 
is surjective, and (ii) follows.

If (i) holds, then the natural map
$\Q_\ell[N] \to \End_{\Q_\ell}(W)$ is surjective
by dimension arguments, and (a) follows.

Part (b) is an immediate corollary of (ii) and
Exercise 15.3 
of \cite{SerreRep}.

The form $f$ in (c) induces an $N$-invariant 
pairing
$\bar{f}: S/\ell S \times S/\ell S \to \F_{\ell}$.
The left and right kernels of $\bar{f}$ are 
$N$-submodules of the simple $N$-module $S/\ell S$,
and thus are $0$, since $\bar{f}\ne 0$.
It follows that $f$ is perfect.
\end{proof}

\begin{defn}
\label{wrded}
We will say a $\Z_\ell[N]$-module $S$ is {\em well-rounded}
if $S$, $N$, and $\ell$ satisfy the hypotheses and the
equivalent conditions (i)-(iii) of Lemma \ref{fulllem}. 
\end{defn}

\begin{rem}
\label{wellprod}
Clearly, if $S_1$ is a well-rounded $\Z_\ell[N_1]$-module 
and $S_2$ is a well-rounded $\Z_\ell[N_2]$-module, then 
$S_1 \otimes_{\Z_\ell} S_2$
is a well-rounded $\Z_\ell[N_1\times N_2]$-module.
\end{rem}

\begin{lem}
\label{lalem}
Suppose $S$ is a well-rounded $\Z_\ell[N]$-module, $W=S\otimes_{\Z_\ell}\Q_\ell$,
and $f$ 
is a non-degenerate $N$-invariant $\Q_\ell$-valued
alternating pairing on $W$.
Suppose $D$ is a group, and
$V$ is a finite-dimensional $\Q_\ell$-vector space and
a simple $D$-module.
Then every ($N \times D$)-stable
$\Z_\ell$-lattice in $W \otimes V$ is of the form
$S \otimes \Gamma$ where $\Gamma$ is a $D$-stable
$\Z_\ell$-lattice in $V$, and
every ($N \times D$)-invariant $\Q_\ell$-valued
pairing on $W \otimes V$ is of the form
$f \otimes h$ where $h$ is a $D$-invariant 
$\Q_\ell$-valued pairing on $W$.
\end{lem}

\begin{proof}
Let $r=\dim(V)$ and $2d=\dim(W)$.
Let $M$ denote the image of
$\Z_\ell[N]\to\End_{\Q_\ell}(W)$.
By Lemma \ref{fulllem}iii,
$M=\End_{\Z_\ell}(S)\cong\M_{2d}(\Z_\ell)$.
Suppose $\Lambda$ is an ($N \times D$)-stable
lattice in $W \otimes V$. Then $\Lambda$ is a
$\Z_\ell[N]$-module and thus an $M$-module.
Since every finitely generated $M$-module is a direct
sum of copies of $\Z_\ell^{2d}$ ($\cong S$) with the
standard action, therefore 
$\Lambda \cong S^r = S \otimes_{\Z_\ell}\Z_\ell^r$,
as $M$-modules.
Since $D$ commutes with $N$, it commutes with
the image of $\Z_\ell[N]$ in $\End(\Lambda)$, i.e.,
with $\End(S)\otimes 1$.
Therefore the image of $D$ in $\Aut(\Lambda)$ is
contained in $1\otimes\Aut(\Z_\ell^r)$, and the
isomorphism $\Lambda \cong S \otimes_{\Z_\ell}
\Z_\ell^r$ is also an isomorphism of $D$-modules. 
Therefore $\Lambda \cong (\Z_\ell^r)^{2d}$ as 
$D$-modules, so $V^{2d} \cong W\otimes V \cong
\Lambda\otimes\Q_\ell \cong (\Q_\ell^r)^{2d}$
as $D$-modules. Since $V$ is a simple $D$-module,
therefore $V^{2d}$ and $(\Q_\ell^r)^{2d}$
are semisimple. It follows that
$V\cong\Q_\ell^r$ as $D$-modules,
and we can identify $\Z_\ell^r$ with a $D$-stable
lattice $\Gamma$ in the $D$-module $V$.

If $g$ is an ($N \times D$)-invariant 
$\Q_\ell$-valued pairing on 
$W \otimes V$, then $W \otimes V\cong V^{2d}$ 
is self-dual with respect to $g$, 
since $V$ is simple as a $D$-module.
Therefore $V$ is self-dual, i.e., admits a
non-degenerate $D$-invariant $\Q_\ell$-valued
pairing $h_0$.
Since $W$ is absolutely simple as an $N$-module,
every $N$-invariant $\Q_\ell$-valued pairing 
on $W$ is a multiple of $f$.
Similarly, since 
$\End_{N\times D}(W\otimes V)=\End_D(V)$,
thus every 
($N \times D$)-invariant $\Q_\ell$-valued
pairing on $W \otimes V$ is of the form
$f \otimes h$, where for some $u\in\Aut_D(V)$ we have
$h(x,y)=h_0(ux,y)$ for every $x, y \in V$.
Since $h_0$ is $D$-invariant, so is $h$.
\end{proof}

The next result extends
Theorem 6.2 of \cite{sympl}. In that result,
$N$ was restricted to being a finite group of 
order not divisible by $\ell$,
$\ell$ was odd, and $D$ was $\bmu_\ell$.

\begin{thm}
\label{repnthythm}
Suppose $S$ is a well-rounded $\Z_\ell[N]$-module,
$W=S\otimes_{\Z_\ell}\Q_\ell$,
and $f$ 
is a non-degenerate $N$-invariant $\Q_\ell$-valued
alternating pairing on $W$.
Suppose $D \subset \GL_m(\Q_\ell)$ 
is an irreducible subgroup.
Assume that there does not exist a $D$-stable 
$\Z_\ell$-lattice $\Gamma$ in $\Q_\ell^m$ with
a perfect $D$-invariant $\Z_\ell$-valued
symmetric pairing.
Suppose $T$ is an ($N \times D$)-stable
$\Z_\ell$-lattice in $W^m$ and
$e$ is a non-degenerate ($N \times D$)-invariant 
$\Z_\ell$-valued alternating pairing on $T$.
Then $e$ is not perfect,
and $\#(T^*/T)=\ell^{2dt}$ for some positive integer
$t$, where $T^*$ is the dual lattice of $T$ with
respect to $e$ and $2d=\dim_{\Q_\ell}(W)$.
\end{thm}

\begin{proof}
Replacing $f$ by a suitable multiple, we may suppose
that $f(S,S)=\Z_\ell$.
By Lemma \ref{lalem}, we have
$T=S \otimes \Gamma$ where $\Gamma$
is a $D$-stable lattice in $\Q_\ell^m$, and
$e=f \otimes h$ where
$h$ is a non-degenerate $D$-invariant
$\Q_\ell$-valued pairing on $\Q_\ell^m$.
Since $e$ is alternating, $h$ is symmetric.
Since $f(S,S)=\Z_\ell$, we have
$h(\Gamma, \Gamma) \subseteq \Z_\ell$.
Since
$h: \Gamma \times \Gamma \to \Z_\ell$
is not perfect, there exists
$z \in \Gamma - \ell \Gamma$
such that $h(z,\Gamma) \subseteq \ell\Z_\ell$.
If $x \in S - \ell S$, then
$x \otimes z \in T - \ell T$ but
$e(x\otimes z,T) \subseteq \ell\Z_\ell.$
Therefore $e$ is
not perfect, so $T^* \ne T$.
The Jordan-H\"older factors of
the $N$-module $T/\ell T \cong (S/\ell S)^m$
are all isomorphic to $S/\ell S$. Since
$T \subset T^* \subseteq \ell^{-r} T$
for some positive integer $r$, the
Jordan-H\"older factors of the $N$-module 
$T^*/T$ are all isomorphic to $S/\ell S$.
Therefore $\#(T^*/T)$ is a power of 
$\#(S/\ell S)= \ell^{2d}$ as desired.
\end{proof}

The symmetric group $S_n$ acts on $\Q_\ell^n$ with
standard basis $\{e_1,\ldots, e_n\}$.
Let
$$U=\{\sum_{i=1}^n c_i e_i \mid \sum_{i=1}^n
c_i=0\},$$
a hyperplane 
in $\Q_\ell^{n}$.
Write $h$ for the restriction to $U$ of the 
standard symmetric pairing on $\Q_\ell^n$. 
Then $h$ is non-degenerate.
There are two well-known examples
of $S_n$-stable $\Z_\ell$-lattices in $U$.
First, the root lattice $Q$ generated
by  ``simple roots''
$\alpha_i=e_i-e_{i+1}$ for $1 \le i \le n-1$.
Second, its dual
with respect to $h$, the weight lattice
$P$ generated by 
$\omega_i=e_i - (e_1+\cdots+e_n)/n$
for $1 \le i \le n$.
We have
$h(\omega_i, \alpha_j)=\delta_{ij}$.
If $t_i \in S_n$ is the
transposition $(i, i+1)$ then
$t_i(x)=x-h(x,\alpha_i)\alpha_i$ for $x \in U$.

\begin{lem}[see \cite{Craig} and Theorem 1 of
\cite{Feit}]
\label{Craiglem}
\begin{enumerate}
\item[(i)] $Q \subseteq P$.
\item[(ii)] $P/Q$ is a cyclic
group whose order is the largest
power of $\ell$ dividing $n$.
\item[(iii)] If $L$ is an
$S_n$-stable $\Z_\ell$-lattice in $U$,
then there exists an integer $t$ such that
  $Q\subseteq \ell^{t}L\subseteq P$.
\item[(iv)] If $n>2$ and $n$ is divisible by $\ell$,
then there are no $S_n$-invariant perfect
$\Z_\ell$-valued
pairings on $P$ or on $Q$.
\item[(v)] If $n>2$, and $\ell$ divides $n$ but
$\ell^2$ does not divide $n$, 
then there does not exist an
$S_n$-stable $\Z_\ell$-lattice in $U$ with a
perfect $S_n$-invariant $\Z_\ell$-valued 
pairing.
\item[(vi)] If $\ell$ does not divide $n$, then 
$Q=P$ and the $\Z_\ell[S_n]$-module $P$ 
is well-rounded.
\end{enumerate}
\end{lem}

\begin{proof}
Since $\alpha_i=\omega_i-\omega_{i+1}$, we have (i).
The group isomorphism $P/Q \cong \Z_\ell/n\Z_\ell$
defined by $\omega_i \mapsto 1$ gives (ii).
For (iii), multiplying $L$ by a suitable
power of $\ell$, we may assume that
$L$ contains $Q$ but does not contain
$\ell^{-1}Q$.
Assume that $L$ is not contained in $P$.
Then there exist $x \in L$ and 
$j\in\{1,\ldots,n-1\}$ such that 
$h(x,\alpha_j)\notin\Z_\ell$.
Multiplying $x$ by a suitable $\ell$-adic integer,
we may assume that $h(x,\alpha_j)=\ell^{-1}$.
Since $x\in L$ and $t_j(x)=x-\ell^{-1}\alpha_j$, 
therefore $\alpha_j/\ell \in L$.
Since the orbit $S_n(\alpha_j)$ contains all the 
roots and $L$ is $S_n$-stable, thus 
$\ell^{-1}Q \subseteq L$.
This contradiction gives (iii).
The $\Q_\ell[S_n]$-module $U$ is absolutely simple, so
every $S_n$-invariant $\Q_\ell$-valued
pairing on $U$ is of the form
$ch(x,y)$ for some $c \in \Q_\ell$.
We have $h(Q,Q) =\Z_\ell$, since
$h(\alpha_1,\alpha_2)=-1$.
Further, $ch(Q,Q)=\Z_\ell$ if and only if $c \in
\Z_\ell^*$.
Suppose that $n$ is divisible by $\ell$.
Since the dual of $Q$ with respect to $h$ is
$P \ne Q$, thus $ch(x,y)$ is never perfect.
Let $\ell^r$ be the largest power of 
$\ell$ dividing $n$.
Since $h(\omega_1, \omega_1)=(n-1)/n$, we have
$\ell^r(P,P)=\Z_\ell$.
Further, $ch(P,P)=\Z_\ell$
if and only if $c \in \ell^r\Z_\ell^*$.
The dual of $P$ with respect to $\ell^rh$
is $\ell^{-r}Q$. If $n>2$, then
$\ell^{-r}Q \ne P$ so $ch(x,y)$ is perfect for no $c$,
and we have (iv).
Part (v) follows from (ii), (iii), and (iv).
Suppose that $n$ is not divisible by $\ell$.
Let bars denote reduction mod $\ell$. 
By (ii), $P=Q$. By (iii), every $S_n$-stable 
$\Z_\ell$-lattice in $U$ is of the form $\ell^t P$.
Therefore the $S_n$-module $\bar{P}=P/\ell P$ 
is simple, by Exercise 15.3 
of \cite{SerreRep}.
It thus suffices to show 
$\End_{S_n}(\bar{P})=\F_\ell$.
Let $\bar{P}_1$ denote the set of elements of
$\bar{P}$ invariant under
$\{\sigma\in S_n \mid \sigma(1)=1\}$.
Then $\bar{\omega}_1$ generates the
one-dimensional $\F_\ell$-vector space
$\bar{P}_1$. 
Suppose $\gamma \in \End_{S_n}(\bar{P})$.
Then $\bar{P}_1$ is $\gamma$-stable, so 
there is a $c \in \F_\ell$ such that 
$\gamma(\bar{\omega}_1)=c\bar{\omega}_1$.
Since $\gamma$ commutes with $S_n$ and the orbit
$S_n(\bar{\omega}_1)$ contains all the
$\bar{\omega}_i$'s, we have
$\gamma(\bar{\omega}_i)=c\bar{\omega}_i$.
Since $\bar{P}$ is generated by the 
$\bar{\omega}_i$'s, we have 
$\gamma=c \in \F_\ell$.
\end{proof}

\section{$\ell$-full abelian varieties}
\label{AVsect}

Suppose $\A$ is an abelian variety over a field $F$
of characteristic $\ne\ell$.

\begin{defn}
\label{ellfull}
We will say that $\A$ is {\em $\ell$-full} 
if the $\Z_\ell[G_{F^{(\ell)}}]$-module 
$T_{\ell}(\A)$ is well-rounded
(see Definition \ref{wrded}).
\end{defn}

\begin{rems}
\label{imcontainssl}
\begin{enumerate}
\item[(i)] 
Since $T_\ell(\A)/\ell T_\ell(\A)=\A_\ell$,
the abelian variety $\A$ is $\ell$-full if and
only if the $G_{F^{(\ell)}}$-module
$\A_\ell$ is absolutely simple.
\item[(ii)] Clearly, $\ell$-fullness is stable under
$F^{(\ell)}$-isogeny. 
\item[(iii)] If $\A$ is $\ell$-full then
$\End_{F^{(\ell)}}(\A)=\Z$.
\item[(iv)] If $\A$ is
$\ell$-full then it has a polarization 
defined over $F$ of degree prime to $\ell$
(this follows from 
Thm.~3 in \S 23 of \cite{Mumford}
and Lemma \ref{fulllem}c).
\item[(v)] Suppose $A$ is $\ell$-full and
$\phi : A \to B$ is an $F$-isogeny
whose degree is a power of $\ell$.
Then $A$ and $B$ are $F$-isomorphic 
(since $\ker(\phi)=A_{\ell^m}$ for some
$m$).
\end{enumerate}
Suppose now that $A$ is an elliptic curve $E$.
\begin{enumerate}
\item[(vi)] 
The image of $G_{F^{(\ell)}}$ in 
$\Aut(T_{\ell}(E))\cong \GL_2(\Z_{\ell})$ 
is the intersection of the image of $G_F$ and 
$\SL(T_{\ell}(E))\cong\SL_2(\Z_{\ell})$. 
In particular, if $G_F \to \Aut(T_{\ell}(E))$ 
is surjective then the image of
$G_{F^{(\ell)}}$ in 
$\Aut(T_{\ell}(E))$ is 
$\SL(T_{\ell}(E))\cong\SL_2(\Z_{\ell})$  
and $E$ is $\ell$-full.
\item[(vii)] 
Suppose $\rho_{E,\ell}(G_{F^{(\ell)}})$ 
contains $\SL_2(\F_\ell)$
(and therefore is $\SL_2(\F_\ell)$).
Then $E$ is $\ell$-full (see (ii) of Lemma
\ref{fulllem}).
\item[(viii)] 
By Lemma 3 on p.~IV-23 of \cite{Serreladicbk},
if $\ell \ge 5$ and $F$ is a number field, then 
the $\ell$-adic representation
$G_{F^{(\ell)}} \to \SL_2(\Z_\ell)$ is
surjective if and only if the mod $\ell$ 
representation 
$\rho_{E,\ell} :
G_{F^{(\ell)}} \to \SL_2(\F_\ell)$
is surjective.
\end{enumerate}
\end{rems}

\begin{ex}
Let $E$ be the elliptic
curve $y^2+y=x^3-x$ over $\Q$. The map
$G_\Q\to\Aut(T_\ell(E))$ is surjective
for all primes $\ell$
(see \cite{SerreInv}, 5.5.6). 
In particular, $E$ is $\ell$-full for all primes
$\ell$.
\end{ex}

\begin{prop}
\label{imcontains}
Suppose $\ell$ is an odd prime, and 
$E$ is an elliptic curve over a field $F$
of characteristic $\ne\ell$.
Suppose $\rho_{E,\ell}(G_F)$ 
contains $\SL_2(\F_\ell)$. Then $E$ is $\ell$-full.
\end{prop}

\begin{proof}
The commutator subgroup $[G_F,G_F]$ lies in 
$G_{F^{(\ell)}}$.
First suppose $\ell>3$.
Then
$\rho_{E,\ell}(G_{F^{(\ell)}})$ 
contains
$[\SL_2(\F_\ell),\SL_2(\F_\ell)]=\SL_2(\F_\ell)$,
and we are done by Remark \ref{imcontainssl}vii.
Suppose $\ell=3$. We have
$[\GL_2(\F_3):\SL_2(\F_3)]=2$.
Let $H=\rho_{E,3}(G_{F^{(3)}})$.
By Remark \ref{imcontainssl}vii,
we may suppose that $H$ is a proper subgroup of
$\SL_2(\F_3)$.
If $\rho_{E,3}(G_F)=\GL_2(\F_3)$, then
$H\supseteq [GL_2(\F_3),\GL_2(\F_3)]=\SL_2(\F_3)$.
Thus, we may assume $\rho_{E,3}(G_F)=\SL_2(\F_3)$.
Then $F$ contains a primitive cube root of unity
and therefore $G_F/G_{F^{(3)}}$ 
is a pro-cyclic $3$-group.
Then $\SL_2(\F_3)/H$ is a cyclic $3$-group.
Since $\#\SL_2(\F_3)=24$,
$H$ has order $8$ which is prime to $3$ and therefore 
the $H$-module $E_3\cong\F_3^2$ is semisimple.
Since $H$ is non-commutative and $\dim_{\F_3}E_3=2$,
it follows that the $H$-module 
$E_3$ is absolutely simple.
\end{proof}

\begin{prop}
\label{ell2}
Suppose $E$ is an elliptic curve
over a number field $F$.  Suppose
$j > 2$ is an integer
such that $\zeta_{2^j} \notin F(\zeta_{2^{j-1}})$, 
and suppose $\rho_{E,2^j}(G_F)$
contains $\SL_2(\Z/{2^j}\Z)$. Then
$\rho_{E,2}(G_{F^{(2)}})=\SL_2(\F_2)$, and
therefore $E$ is $2$-full.
\end{prop}

\begin{proof}
Write $F'=F(\zeta_{2^{j-1}})$ and 
$F''=F(\zeta_{2^j})$.
The mod $2$ representation 
$$\rho_{E,2} : G_F \to \Aut(E_2) \cong \SL_2(\F_2) 
\cong S_3$$ is surjective, 
since $\SL_2(\Z/{2^j}\Z) \to\SL_2(\F_2)$ is
surjective. 
Let $G=\rho_{E,2}(G_{F^{(2)}})$.
Since $G_F/G_{F^{(2)}}$ is a pro-$2$-group,
the index of $G$ in $\SL_2(\F_2) \cong S_3$ 
is a power of $2$. Thus either $G\cong\SL_2(\F_2)$ 
and we are done, or $G\cong A_3$.
Assume the latter. We have
\begin{equation}
\label{fjeqn}
\rho_{E,2^j}(G_{F'})\supseteq\rho_{E,2^j}(G_{F''})=
\rho_{E,2^j}(G_{F})\cap\SL_2(\Z/2^j\Z)=\SL_2(\Z/2^j\Z).
\end{equation}
Let $\sigma$ be the composition
$$G_{F'} \stackrel{\rho_{E,2}}{\twoheadrightarrow}
\GL_2(\F_2)=S_3 \to S_3/A_3=\{1,-1\},$$
and let $F_\sigma$ be the fixed field of
$\ker(\sigma)$.
Then $F_\sigma$ is a quadratic extension of $F'$.
By (\ref{fjeqn}), $F_\sigma\ne F''$. 
Since $j>2$, $F''$ is the only quadratic extension of
$F'$ in $F^{(2)}$. Therefore $F_\sigma$ is not
contained in $F^{(2)}$, so 
$\sigma(G_{F^{(2)}})\ne 1$.
This contradicts the assumption that $G=A_3$.   
\end{proof}

\begin{thm}
\label{modN}
Suppose $F$ is a number field and $n$ is a positive
integer.
Then there are infinitely many elliptic curves $E$
over $F$, non-isomorphic over $\C$, 
for which $\SL_2(\Z/n\Z)\subseteq\rho_{E,n}(G_F)$.
\end{thm}

\begin{proof}
As shown on pp.~145--6 of \cite{SerreMW}, there is an
elliptic curve $E(t)$ over the function field $\Q(t)$,
with $j$-invariant $t$, such that
$$\Gal(\Q(t,E(t)_n)/\Q(\zeta_n)(t))=\SL_2(\Z/n\Z)$$
and $\Q(\zeta_n)$ is algebraically
closed in $\Q(t,E(t)_n)$.
Thus 
$$F(\zeta_n)(t)\cap \Q(t,E(t)_n) = \Q(\zeta_n)(t),$$
and we have 
$$\Gal(F(t,E(t)_n)/F(\zeta_n)(t))=\SL_2(\Z/n\Z).$$
By Prop.~2 on p.~123 of \cite{SerreMW}
and Hilbert's Irreducibility Theorem 
(p.~130 of \cite{SerreMW}), there are infinitely many 
specializations $t_0\in F$ such that 
$$\Gal(F(E(t_0)_n)/F(\zeta_n))=
\Gal(F(t,E(t)_n)/F(\zeta_n)(t))=\SL_2(\Z/n\Z),$$
and therefore 
$\Gal(F(E(t_0)_n)/F)\supseteq \SL_2(\Z/n\Z)$.
\end{proof}

\begin{cor}
\label{manyell}
Let $F$ be a number field and $r$ a positive integer. 
Then there are infinitely many elliptic curves 
over $F$, non-isomorphic over $\C$,
that are $\ell$-full for all prime divisors 
$\ell$ of $r$.
\end{cor}

\begin{proof}
Choose $j>2$ such that $\zeta_{2^j}\notin
F(\zeta_{2^{j-1}})$
(such a $j$ exists since $[F(\zeta_{2^{k}}):F]$ is
unbounded as $k \to \infty$).
Let $n'$ be the product of the prime divisors of $r$.
Let  $n=n'$ if $r$ is odd and let $n=2^{j}n'$ if
$r$ is even.
By Theorem \ref{modN}, there are infinitely many 
elliptic curves $E$ over $F$, non-isomorphic over 
$\C$, for which 
$\SL_2(\Z/n\Z)\subseteq\rho_{E,n}(G_F)$.
Thus for all odd prime divisors $\ell$ of $r$
we have $\SL_2(\F_{\ell})\subseteq\rho_{E,\ell}(G_F)$,
so $E$ is $\ell$-full 
by Proposition \ref{imcontains}.
If $r$ is even, then $E$ is $2$-full
by Proposition \ref{ell2}.
\end{proof}

\begin{prop}[see \cite{SerreMW}]
\label{frat1}
Suppose $Y$ is an abelian variety over a field $F$ of
characteristic $\ne \ell$. 
Then there exists a finite Galois extension 
$K_\ell/F$ such that every finite extension of 
$F$ linearly disjoint from $K_\ell$ is linearly
disjoint from $F_Y^{(\ell)}$. 
\end{prop}

\begin{proof} 
Let $H=\Gal(F_Y^{(\ell)}/F)$, 
an $\ell$-adic Lie group.
By the Proposition 
and Example 1 in \S~10.6 
of \cite{SerreMW}, there exists an open normal 
subgroup $\Phi\subset H$ of finite index
such that whenever $\tilde{H}$ is a closed normal 
subgroup in $H$ with $H=\tilde{H}\Phi$, 
then $\tilde{H}=H$.
Let  $K_\ell$ be the fixed field of $\Phi$. 
Suppose $K$ is a finite extension of $F$, 
let $\tilde{H}=\Gal(K_Y^{(\ell)}/K)$, and
view $\tilde{H}$ as a closed normal 
subgroup of $H$.
If $K$ and $K_\ell$ are linearly disjoint,
then the restriction to $\tilde{H}$ of
$H\to H/\Phi$
is surjective, i.e., $H=\tilde{H}\Phi$.
Therefore $\tilde{H}=H$, so $K$ and 
$F_Y^{(\ell)}$ are linearly
disjoint.
\end{proof}

\section{A general method}

Suppose $\A$ is an abelian variety over a field
$F$ of characteristic $\ne \ell$.
There is a canonical $\Z_{\ell}$-bilinear
$G_F$-equivariant perfect pairing 
$$e_{\ell}: T_{\ell}(\A) \times T_{\ell}(\A^*) \to
\Z_{\ell}(1)\cong\Z_{\ell}$$
where $\Z_\ell(1)$ is the projective limit of the
groups of $\ell^m$-th roots of unity
(\cite{Mumford}, \S 20; \cite{Lang}, Ch.\ VII, \S 2). 
Since $G_{F^{(\ell)}}$ acts trivially on
$\Z_\ell(1)$, the pairing $e_{\ell}$ is
$G_{F^{(\ell)}}$-invariant.
There exists an $F$-isogeny $A \to A^*$, and therefore 
a surjection 
$A_{\ell^\infty}\twoheadrightarrow A^*_{\ell^\infty}$
defined over $F$. Thus, 
$F_{A^*}^{(\ell)}\subseteq F_{A}^{(\ell)}$.
Since $e_{\ell}$ is perfect, it follows that
$$F^{(\ell)} \subseteq F_{A}^{(\ell)}.$$

Suppose $\mathcal{L}$ is an invertible sheaf on
$\bar{\A}$.
Let $\chi(\mathcal{L})$ denote the Euler
characteristic, and let 
$\phi_\mathcal{L} : \bar{\A} \to \bar{\A^*}$ be
the associated natural homomorphism
defined in \S 13 (see p.~131) of \cite{Mumford}.
Then $\deg(\phi_\mathcal{L})=\chi(\mathcal{L})^2$ 
(see \S 16 of \cite{Mumford}), 
and $\phi_\mathcal{L}$ is an isogeny
if and only if $\chi(\mathcal{L}) \ne 0$.
If $\mathcal{L}$ is ample, then 
$\chi(\mathcal{L}) > 0$ and the isogeny 
$\phi_\mathcal{L}$ is called a polarization on
$\bar{\A}$.

Suppose now that $\phi_\mathcal{L}$ is defined over
$F$ (i.e., is obtained by extension of scalars from 
a morphism $\A \to \A^*$).
Then $\phi_\mathcal{L}$ induces a homomorphism of
$G_F$-modules
$f_\mathcal{L}: T_\ell(\A) \to T_\ell(\A^*)$
which is injective if and only if $\phi_\mathcal{L}$
is an isogeny. If $f_\mathcal{L}$ is injective, 
then the index of
the image of $T_\ell(\A)$ in $T_\ell(\A^*)$ 
is the exact power of $\ell$ dividing 
$\deg(\phi_\mathcal{L})$ 
(\cite{Lang}, Ch.\ VII, \S 1, Thm.\ 1 and its proof).
Extending $\Q_\ell$-linearly, the induced homomorphism 
$f_\mathcal{L}: V_\ell(\A) \to V_\ell(\A^*)$
is an isomorphism of $\Q_{\ell}[G_F]$-modules.
The homomorphism $f_\mathcal{L}$ gives rise to an
alternating $G_F$-equivariant pairing
$$E^\mathcal{L}: T_{\ell}(\A) \times T_{\ell}(\A) \to
\Z_{\ell}(1)\cong\Z_{\ell}, \quad x,y \mapsto
e_{\ell}(x,f_\mathcal{L}(y))$$
(see \S 20 of \cite{Mumford}). 
The form $E^\mathcal{L}$ is non-degenerate 
if and only if
$\phi_\mathcal{L}$ is an isogeny. 
It is perfect if and only if $\deg(\phi_\mathcal{L})$ 
is not divisible by $\ell$. 

Assume that $\phi_\mathcal{L}$ is an isogeny.
Let $T_\ell(\A)^*$ denote  
the dual lattice of  $T_\ell(\A)$ with
respect to $E^\mathcal{L}$. Then 
$$T_{\ell}(\A)\subseteq
T_\ell(\A)^*=f_\mathcal{L}^{-1}(T_{\ell}(\A^*))
\subset V_{\ell}(\A),$$
and $\#(T_\ell(\A)^*/T_\ell(\A))=
\#(T_{\ell}(\A^*)/f_\mathcal{L}(T_{\ell}(\A))$
is the exact power of $\ell$ dividing
$\deg(\phi_\mathcal{L})$
(see the proof of Thm.\ 3 in \S 20 of
\cite{Mumford}). 

\begin{lem}
\label{twistgalois}
Suppose $Y$ is an abelian variety over a field $F$
of characteristic $\ne \ell$, and 
$s$ is a positive integer. 
Suppose $K/F$ is a finite Galois extension 
with $\Gal(K/F)=D\subset\GL_s(\Z)$, such that
$-1\notin D$.  
Let $\A$ be the $K/F$-form of $Y^s$ 
attached to the inclusion
$\Gal(K/F)=D\subset\GL_s(\Z)\subseteq \Aut(Y^s)$. 
Let $M$ denote the image of $G_F$ 
in $\Aut(T_\ell(Y))$. Then:
\begin{enumerate}
\item[(a)] the natural map 
$f : M \times D \to 
\Aut(T_\ell(Y)\otimes_{\Z_\ell}\Z_\ell^s)
=\Aut(T_\ell(\A))$
is injective, 
\item[(b)] $F_A^{(\ell)}=KF_Y^{(\ell)}$.
\end{enumerate} 
Suppose now that $K$ and $F_Y^{(\ell)}$ 
are linearly disjoint over $F$. Then:
\begin{enumerate}
\item[(c)] the image of $G_F$  
in $\Aut(T_\ell(A))$ is $f(M\times D)$, 
\item[(d)] if $Y$ is $\ell$-full and the
$\Z_\ell[D]$-module $\Z_\ell^s$ is well-rounded, 
then $\A$ is $\ell$-full.
\end{enumerate}
\end{lem}

\begin{proof}
The kernel of
$$\Aut(T_\ell(Y))\times \GL_s(\Z_\ell)\to
\Aut(T_\ell(Y)\otimes_{\Z_\ell}\Z_\ell^s)
=\GL_s(\End_{\Z_\ell}(T_\ell(Y))$$
is $\{(a, a^{-1})\mid a \in \Z_\ell^*\}$.
Since $-1\notin D$, we have
$$1=D\cap\Aut(T_\ell(Y))
\subseteq \GL_s(\End_{\Z_\ell}(T_\ell(Y)),$$
and (a) follows. From the natural injections 
$$\Gal(KF_Y^{(\ell)}/F) \hookrightarrow
M \times D \stackrel{f}{\hookrightarrow}
\Aut(T_\ell(A))$$
we have (b).
Suppose $K$ and $F_Y^{(\ell)}$ 
are linearly disjoint over $F$. Then
the natural injection
$\Gal(KF_Y^{(\ell)}/F) \hookrightarrow
\Gal(F_Y^{(\ell)}/F)\times 
\Gal(K/F) \cong M \times D$ 
is surjective, and (c) follows.
Now (d) follows from Remark \ref{wellprod}, 
by applying (c) and (a) with $F^{(\ell)}$ in 
place of $F$. 
\end{proof}

\begin{thm}
\label{avprop}
Suppose $F$ is a field
of characteristic different from $\ell$,
$Y$ is an $\ell$-full $d$-dimensional
abelian variety over $F$,
and $s$ is a positive integer. Suppose
$D$ is an irreducible subgroup of $\GL_s(\Z)$
such that $-1\notin D$ and such that
there does not exist a $D$-stable
$\Z_\ell$-lattice in $\Q_\ell^s$
with a perfect $D$-invariant $\Q_\ell$-valued 
symmetric pairing.
Suppose $K/F$ is a Galois extension with $\Gal(K/F)=D$
and with $K$ and $F_Y^{(\ell)}$ linearly
disjoint over $F$.
Let $A$ be the twist of $Y^s$ via
$$\Gal(K/F)=D \subset \GL_s(\Z) \subseteq \Aut(Y^s),$$
and suppose $B$ is an abelian variety
$F^{(\ell)}$-isogenous to $A$.
Then  $\ell^{2d}$ divides the degree of every
$F^{(\ell)}$-polarization on $B$,
and $\ell^d$ divides the Euler characteristic of every
invertible sheaf $\mathcal{L}$ on $\bar{B}$ such that
$\phi_\mathcal{L}$ is defined over $F^{(\ell)}$.
If $\chi(\mathcal{L})\ne 0$ and $\ell^r$ 
is the highest power of
$\ell$ dividing $\chi(\mathcal{L})$, then
$\ell^r$ is a power of $\ell^{d}$.
\end{thm}

\begin{proof}
If $\phi_{\mathcal{L}}$  is not an
isogeny then $\chi(\mathcal{L})=0$ 
and we are done.
Suppose now that $\phi_\mathcal{L}$  is an isogeny.
Let $N$ denote the image of $G_{F^{(\ell)}}$ 
in $\Aut (T_\ell(Y))$.
By Lemma \ref{twistgalois},
the image of $G_{F^{(\ell)}}$ in
$\Aut(V_\ell(A))$ is $N \times D$.
Therefore the $G_{F^{(\ell)}}$-stable 
$\Z_\ell$-lattice $T_\ell(B)$
is also ($N \times D$)-stable.
Since $E^\mathcal{L}$ is 
$G_{F^{(\ell)}}$-invariant, 
it is also ($N \times D$)-invariant.
Applying Theorem \ref{repnthythm} with
$S=T_\ell(Y)$, $T=T_\ell(B)$, and
$e=E^\mathcal{L}$,
then $\#(T^*/T)=\ell^{2dt}$ for some 
positive integer $t$, and the result follows.
\end{proof}

\section{Application}
\label{thm1sect}

\begin{thm}
\label{eulerthm}
Suppose $F$ is a number field and $\r$ is a positive 
integer. Then there exists an abelian variety $A$ 
over $F$ with the following properties:
\begin{enumerate}
\item[(a)] if $B$ is an abelian variety $F$-isogenous
to $A$ and $\mathcal{L}$ is an invertible sheaf 
on $\bar{B}$ such that $\phi_\mathcal{L}$ is defined 
over $F$, then $\chi(\mathcal{L})$ is divisible by $\r$;
\item[(b)] if $\ell$ is a prime, 
$\ell^{a}$ divides $\r$, 
$B$ is an abelian variety $F^{(\ell)}$-isogenous
to $A$, and ${\mathcal{L}}$ is an invertible sheaf 
on $\bar{B}$ such that $\phi_\mathcal{L}$ is defined 
over $F^{(\ell)}$, then $\chi(\mathcal{L})$ is 
divisible by $\ell^{a}$.
\end{enumerate}
\end{thm}

\begin{proof}
We may assume that $\r>1$.
Let $m$ be a positive integer
relatively prime to $\r$ and such that for every prime
divisor $\ell$ of $\r$,
$\ell^{m-1}$ does not divide $\r$.
Let $\n=6$ if $\r$ is a power of $2$, and otherwise
let $\n$ be the product of the prime divisors of $\r$.
Identify the symmetric group $S_t$ in the standard 
way as an irreducible subgroup of $\GL_{t-1}(\Z)$. 
By Lemma \ref{Craiglem}, 
for all prime divisors $\ell$ of $r$
the $\Z_\ell[S_m]$-module
$\Z_\ell^{m-1}$ is well-rounded.
Note that $m>2$ and $n>2$, so for $t=m$ and $n$ we have
$-1\notin S_t \subset GL_{t-1}(\Z)$.
By Lemma \ref{Craiglem}v, no $S_\n$-stable 
$\Z_\ell$-lattice in $U\cong\Q_\ell^{n-1}$
has a perfect $S_n$-invariant $\Z_\ell$-valued
pairing.
By Corollary \ref{manyell}, there exists an elliptic
curve $E$ over $F$  that is $\ell$-full for all
prime divisors $\ell$ of $r$.
By Prop.~\ref{frat1}, 
for every prime $\ell$ there exists a finite
Galois extension $K_\ell/F$ such that
every finite extension of $F$ linearly
disjoint from $K_\ell$ is linearly
disjoint from $F_E^{(\ell)}$.
There exists a finite Galois extension $L$ of $F$
with $\Gal(L/F)=S_m\times S_n$
which is linearly disjoint from $K_\ell$ over 
$F$ for all prime divisors $\ell$ of $r$ (see
\cite{SerreGT}, especially Prop.~3.3.3).
We can write $L=K'K$ with $\Gal(K'/F)=S_m$
and $\Gal(K/F)=S_n$. Then $K'$ and $K$ 
are linearly disjoint over $F$,
and $K'$ (resp., $K$) and $F_E^{(\ell)}$ 
are linearly disjoint over $F$. 
Let $Y$ be the $K'/F$-form of $E^{m-1}$ attached to
$$\Gal(K'/F)= S_m \subset 
\GL_{m-1}(\Z) \subseteq \Aut(E^{m-1}),$$
and let $A$ be the $K/F$-form of $Y^{\n-1}$ 
attached to
$$\Gal(K/F)=S_\n \subset \GL_{\n-1}(\Z) 
\subseteq \Aut(Y^{\n-1}).$$
By Lemma \ref{twistgalois}b, 
$F_Y^{(\ell)}=K'F_E^{(\ell)}$, so $K$ and
$F_Y^{(\ell)}$ are linearly disjoint over $F$.
By Lemma \ref{twistgalois}d, $Y$ is $\ell$-full. 
Part (b) now follows from Theorem \ref{avprop} 
with $D=S_n$. 
Part (a) follows from (b) for 
all prime divisors $\ell$ of $r$.
\end{proof}

Applying Theorem \ref{eulerthm} when $\mathcal{L}$
is ample yields Theorem \ref{thm1}.

\end{document}